\documentclass[12pt,reqno]{amsart}
\pdfoutput=1
\usepackage{amsmath}%
\usepackage{amsfonts}%
\usepackage{amssymb}%
\usepackage[margin=.5in]{geometry}
\usepackage{subfig}
\usepackage{url} 
\usepackage{ifpdf}
\ifpdf
\usepackage{graphicx}
\else
\usepackage[draft]{graphicx}
\fi
\usepackage{color}
\usepackage{enumerate}
\usepackage[T1]{fontenc}
\usepackage{caption}
\usepackage{booktabs,multirow}
\usepackage{siunitx}
\usepackage[flushleft]{threeparttable}
\usepackage{enumerate}
\begin{document}

\title[Puiseux Series]{An algorithm for determining the radii of convergence of algebraic power series}
\author{Dominic C. Milioto}
\email{icorone@hotmail.com}
\date{\today}
\begin{abstract}
This paper describes an algorithm for determining radii of convergence of power expansions for algebraic functions and the testing done to check it.  Since the current methods for computing these series are iterative, standard methods for computing radii of convergence cannot in general, be used.  However, relying on geometric properties of algebraic functions, convergence radii of these series can be determined precisely.
\end{abstract}
\maketitle
\section{Introduction}

The objects studied in this paper are fractional power expansions of algebraic functions
\begin{equation}
f(z,w)=a_0(z)+a_1(z)w+a_2(z)w^2+\cdots+a_{n}(z)w^{n}=0,
\label{eqn001}
\end{equation}
such that $f(z,w)$ is irreducible over the rationals with $z$ and $w$ complex variables and the coefficients, $a_i(z)$, polynomials in $z$ with rational coefficients.  The degree of the function is the highest power of $w$.   By the Implicit Function Theorem, this equation defines locally, an analytic function $w(z)$ when $\displaystyle \frac{\partial f}{\partial w}\neq 0$.   And by Newton-Puiseux's Theorem\cite{Nawak}, $(\ref{eqn001})$ can be factored over the field of fractional power series, $K\left(\left\{z\right\}\right)$ as
\begin{equation}
f(z,w)=\prod_{j=1}^{S_r} \left(w-w_d(z)_j\right),\quad w_d(z)\in K\left(\left\{z\right\}\right),
\label{eqn001b}
\end{equation}
where each $w_d(z)$ is a branch, possibly multivalued, of the function.  These fractional power series have a radius of convergence at least equal to the distance to the nearest singular point. However, one result of this paper is to demonstrate the actual radius of convergence can and often does extend beyond many singular points.  In the test cases described below, one branch is shown to extend across $118$ singular points.

Since the current methods for computing fractional power series of algebraic functions rely on iterative means \cite{Poteaux}, radii of convergence, in general, cannot be calculated by the ordinary techniques such as the ratio test or root test.   This paper describes a method to compute these values precisely.   

The purpose of this paper is four-fold:
\begin{enumerate}
\item{Describe in detail, the Newton polygon algorithm as described by Kung \cite{Kung},}
\item {Demonstrate that a numeric version of Newton Polygon can be implemented with satisfactory results,}
\item {Implement an algorithm for determining the radius of convergence of an algebraic power series,}
\item {Provide a software tool for further investigating algebraic functions.}
\end{enumerate}
\section{Some properties of algebraic functions used in this paper}
Fractional power expansions of algebraic functions are called Puiseux series and are in this paper computed by the method of Newton polygon.  These series often represents only a small portion of the function near the origin.  However, the entire function can be represented by Puiseux series in annular discs surrounding the origin.  These annular expansions can in principle, be computed by numeric means using a modified version of Laurent's expansion theorem and if computed accurately, can be used in the algorithms below to determine their annular domain of convergence.  However, this paper focuses only on power series inside the disc $D(0,|r_c|)$.  

The resultant of $f(z,w)$ with it's derivative $f_w$ is denoted by $R(f,f_w)$.  A point where one or both of $a_n(z)$ or $R(f,f_w)$ is zero is a singular point of $f$.   A point where $a_n(z)=0$ is also a pole, possibly ramified, of the function.  In this paper, singular points are labeled as $s_n$.  A property of algebraic functions which distinguishes them from single-valued functions is algebraic functions can be both singular and analytic at a singular point.  That is, a singular point may not affect all coverings of an algebraic function unless the function is fully-ramified at the singular point.  For example, a $10$-degree function may have only a single $2$-cycle branch and eight $1$-cycle coverings at a singular point.  In this case, the $2$-cycle covering is singular.  The eight single-cycle coverings are not analytically affected at this singular point unless one is affected by a pole of the function.  However, if the function fully-ramifies into a $10$-cycle branch at this singular point, all coverings would be affected.  It is for this reason algebraic power series can have a radius of convergence extending beyond the first singular point:  Their branch coverings may simply not be singular at a singular point.  Only when the covering becomes singular does the convergence radius of its power expansion become established. The main objective of this paper is to identify which singular point is interrupting the analyticity of branch cycles thereby establishing the radius of convergence of their power expansions.

 The following conventions are used in this paper:

\begin{enumerate}

\item The Puiseux expansions of the roots of $(\ref{eqn001})$ consist of a set of $d$-valued branches.  A branch is sometimes called a $d$-cycle where $d$ is a positive integer.  A power series in $z^{1/3}$ would be a $3$-cycle branch.  It has three coverings over the complex $z$-plane just like the function $f(z)=\sqrt[3]{z}$.  The sum of the cycles is always equal to the degree of the function in $w$.
 
\item  In order to identify a particular branch, $w(z)$ of $f$, a two-level identification is used.  The first level is the cycle number representing the number of coverings.  The second identifier is the sort order of the branch: A base point $z_m$ on the real axis is selected mid-way between the origin and a disc of radius equal to the distance to the nearest non-zero singular point.  When the function has multiple $d$-cycle branches, the values $\{w_{d_i}(z_m)\}$ are computed and ordered by increasing real part then imaginary part.  The first value in the sort order for each branch is identified and the branches are then ordered according to this sort order as $w_{d,1},w_{d,2}$ and so forth.
  
\item The concept of ``branch'' is used throughout this paper and differs in meaning to the definition encountered in basic Complex Analysis.  In this paper, the term ``branch'' refers to a multi-valued $d$-cycle root of $f(z,w)$ given by $(\ref{eqn001b})$. 

\item The following discussion makes use of the term, ``extending a branch over a singular point''.  This is in reference to the discussion above about singularities and the coverings they affect.

\item In this paper, $r_c$ is a positive integer representing a singular ring when non-zero singular points are arranged in order of increasing absolute value.  The smallest non-zero singular point is therefore on ring one, singular points with the next largest absolute value are on ring two and so forth.  $|r_c|$ however, represents the absolute value of the associated singular points.

\end{enumerate}

\section{Newton polygon method}

Given $(\ref{eqn001})$, we seek to compute power expansions
\begin{equation}
w_d(z)=\frac{1}{z^E}\sum_{n=0}^{\infty}c_n \left(z^{1/d}\right)^n
\label{eqn002}
\end{equation}  
where $E$ is the normal exponent and $d$ is the cycle number, for all branches of the function. The simplest method of computing $(\ref{eqn002})$ is by the method of Newton polygon.  For a function of degree $n$ in $w$, the method computes a set of $n$ series representing the roots of $(\ref{eqn001})$.   However, the actual branches of the function are often multivalued and ramify or wind around the origin multiple times.  For example, a ten-degree function can ramify into a single-valued branch, a $2$-cycle branch, a $3$-cycle, and $4$-cycle branch.  Four of the power series computed by the method would be conjugates of the $4$-cycle, three series, conjugates of the $3$-cycle, and so forth.   Thus in the case of the $4$-cycle branch, the algorithm produces the set
\begin{equation}
\{w_{4}(z)_j\}=\frac{1}{z^E}\sum_{n=0}^{\infty} c_n\left(e^{2j\pi i/4}\right)^n \left(z^{1/4}\right)^n,\quad j=0,1,2,3.
\label{eqn0010}
\end{equation}
In the particular implementation of the algorithm presented here, a basis set of power series is produced consisting of a single power series for each $d$-cycle branch.  The remaining conjugate series can be computed from the basis series using the appropriate form of  (\ref{eqn0010}).

\begin{figure}
	\centering
		\includegraphics[scale=0.75]{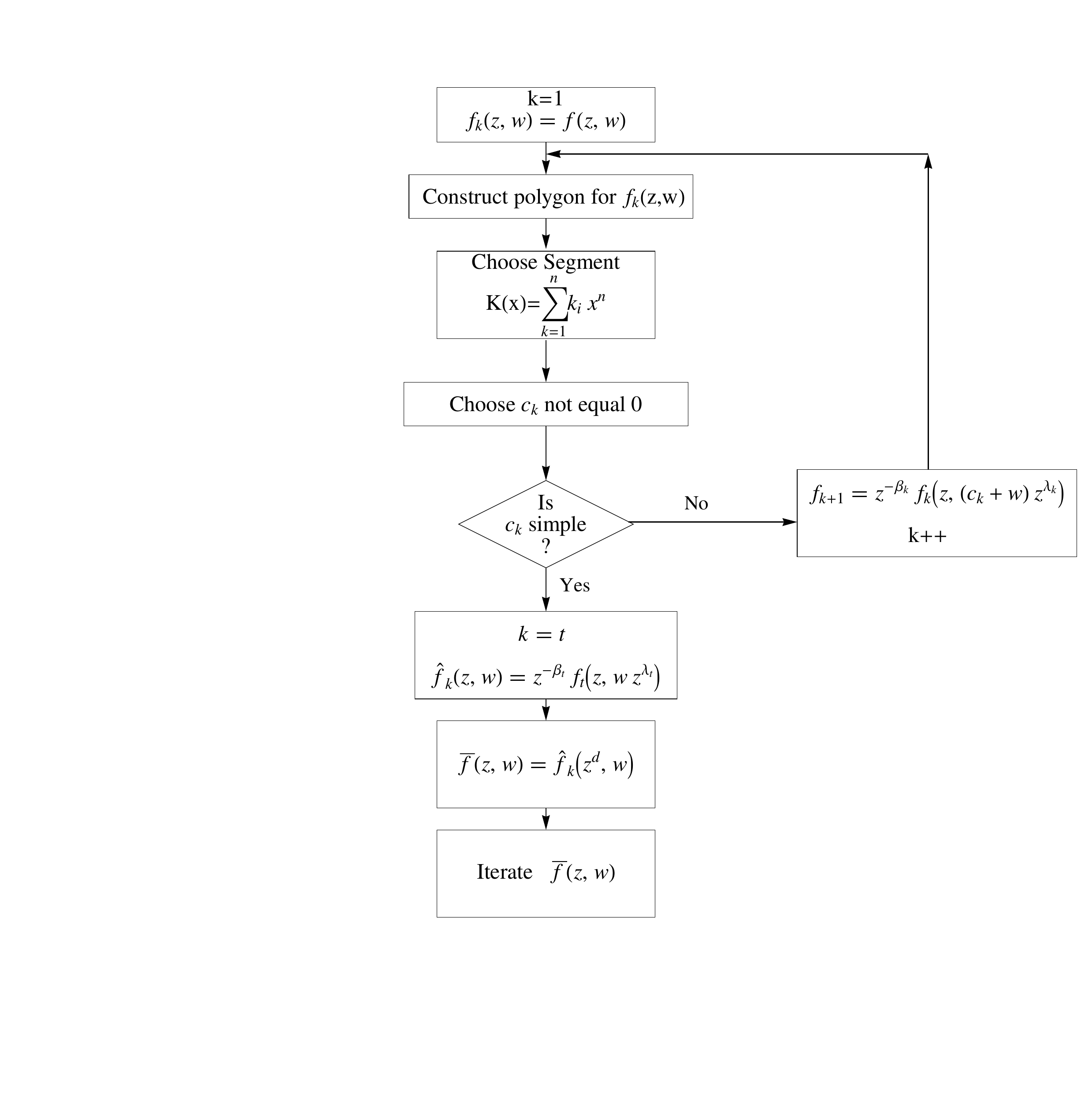}
		\caption{Flowchart for Newton Polygon}
	\label{fig0010}
\end{figure}

A theoretical basis for the Newton polygon algorithm can be found in Walker \cite{Walker}.  The particular implementation of the algorithm for this paper is taken from Kung \cite{Kung}, a flowchart of which is shown in Figure $\ref{fig0010}$.  We begin first by representing a branch as
\begin{equation}
w(z)=c_1 z^{\lambda_1}+c_2 z^{\lambda_1+\lambda_2}+c_3 z^{\lambda_1+\lambda_2+\lambda_3}+\cdots
\label{eqn038}
\end{equation}
in which $\lambda_1\geq 0$ and $\lambda_i>0$ for $i>1$ and $c_i\neq 0$.  The steps in the process are as follows:
\begin{enumerate}
\item{Normalizing the function}:

In order to achieve the requirement above on the exponents $\lambda_i$, the function is normalized to eliminate negative exponents in the expansion. The normalization process removes algebraic poles at the origin from the function which arise when zero is a root of $a_n(z)$.  For example,
$$
f(z,w)=(78+49 z-46 z^2)\text{}+(-2 z)w+(91-34 z-80 z^2)w^2+(52 z+47 z^2)w^3
$$
has a pole at the origin and thus one branch of the function would have a negative exponent in its expansion.  The normalization process produces for this function, 
$$
\begin{aligned}
g(z,w)&=z^2 f(z,w/z)\\
&=(78 z^2+49 z^3-46 z^4)\text{}+(-2 z^2)w+(91-34 z-80 z^2)w^2+(52+47 z)w^3,
\end{aligned}
$$
which in this case $E=1$.  We can then expand $g(z,w)$ around the origin, which now has no singular point there and hence no expansions with negative exponents, and then multiply each series by the quantity $1/z$ to arrive at the expansions for $f(z,w)$.
\item{Computing the Newton polygon and extracting the lower Newton leg:}

The principle device for beginning the computation of the coefficients $c_n$ for the series is the Newton polygon which for the function $f(z,w)$ is the convex hull of it's support.  The support of the function is a set of points, $\{i,\alpha(a_i)\}$  where $\alpha(a_i)$ is the order of $a_i(z)$.  The order of $a_i(z)$ is the lowest power of $z$ in $a_i(z)$ with the exception that if $a_i(z)=0$, then $\alpha(a_i)=\infty$ and the point $(i,\infty)$ is omitted.  Drawing the convex hull of these points and taking the left-most lowest segments so that all remaining points are either above or to the right of the segments then gives the lower Newton leg.  The Newton polygon for 
\begin{equation}
\begin{aligned}
f(z,w)&=(z^{14}+3 z^{15})\text{}+(2 z^{15}+2 z^{16})w+(3 z^{15}-20 z^{16})w^2\\
&+(4 z^{15})w^3+(10 z^5-z^6+2 z^7)w^4+(8 z^{10})w^5+(9 z^{10})w^6+(20 z)w^7\\
&+(3 z)w^8+(2)w^9+(5)w^{10}
\end{aligned}
\label{eqn039}
\end{equation}
is shown in Figure $\ref{fig:polygon1}$.  In this example, the lower Newton leg has four segments given by the red lines.  The blue dashed lines provide the slope and vertical intercept for each segment which gives the exponents $\beta$ and $\lambda$ discussed below.
\begin{figure}
	\centering
		\includegraphics{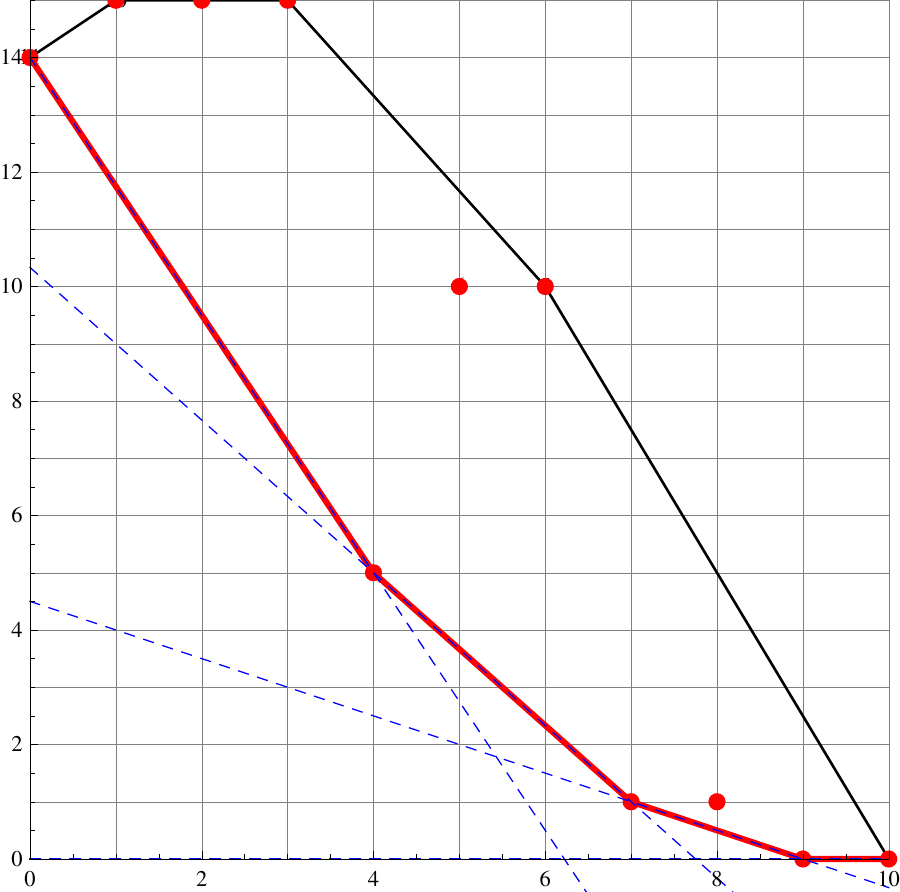}
	\caption{Newton polygon for $(\ref{eqn039})$}
	\label{fig:polygon1}
\end{figure}
\item{Deriving the characteristic equation:}

Choose a segment from the lower Newton leg and form the characteristic equation of the segment.  The principle is based on letting $w(z)=z^{\lambda_1}\left(c_1+w_1\right)$ and then making the substitution $f(z,z^{\lambda_1}(c_1+w_1))$.  Extracting the lowest power of $z$ alone and setting it to zero then provides a means of computing $c_1$. However, the form of the characteristic equation can be read off directly from the polygon diagram:  for each segment, choose the points intersected by the segment.  In the case of Figure $\ref{fig:polygon1}$, the top most segment intersects the points $(0,14)$ and $(4,5)$.  These are from the order-14 term of $a_0(z)$ and the order-5 term of $a_4(z)$.  Now choose the coefficient of each of those terms and label them $b_{\text{degree},\text{order}}=b_{i,\alpha}$.  In the two cases above, we have $1_{0,14}$ and $10_{4,5}$.  We can write this set of points as:
$$
\{b_{i_1,\alpha_1},b_{i_2,\alpha_2},\cdots,b_{n_Q,\alpha_Q}\}
$$
Now form the characteristic equation
\begin{equation}
K(x)=\sum_{j=1}^Q b_{i_j,\alpha_j}x^{i_j}=0
\label{eqn033}
\end{equation}
with $Q$ being the total number of points.  This give for the segment above, $\displaystyle K(x)=1+10x^4$.  The $k$ distinct, non-zero solutions to the characteristic equation are the first coefficients $c_1$ of $k$ power series for the function.  If the zeros are multiple roots, then we create a series for each set of multiple roots by setting $c_1$ for each root, and let $f_1(z,w_1)=z^{-\beta_1}f(z,z^{\lambda_1}(c_1+w_1))$ where $\lambda_1$ is the negative of the slope for the segment and $\beta_1$, the vertical intercept of the segment.  We then compute a second-level polygon for $f_1(z,w_1)$ with the exception that segments with zero slopes are now omitted from the Newton leg, derive a second characteristic equation for $c_2$, and continue in this way if necessary for $c_3$ and so forth until we reach $f_k(z,w_k)$ with simple roots for its characteristic equation.
\item{Convert to regular form:}

Once simple roots are obtained for the characteristic equation, we can use a variation of Newton iteration to compute additional terms of the series.  In order to use this method, we need a polynomial with integer powers; the procedure above will often produce an expression with fractional powers.  But we can convert the expression $f_k(z,w_k)$ above to one with integer powers with the following two substitutions:
\begin{equation}
\begin{aligned}
\tilde{f}(z,w)&=\frac{1}{z^{\lambda_k}} f_k(z,z^{\lambda_k} w) \\
\overline{f}(z,w)&=\tilde{f}(z^d,w)
\end{aligned}
\label{eqn034}
\end{equation}
where $k$ is the last recursion of the function which did not produce a characteristic equation with multiple roots, and $d$ is the lowest common denominator of the exponents $\{\lambda_j\},\quad j=1,2,\cdots,k$. 
\item{Normal Iteration:}

The terms in $\overline{f}(z,w)$ not containing a factor of $z$ make up the characteristic equation in $w$ so that the first simple root computed  above, $r_k$, is a zero to $\overline{f}(0,w)$.  That is, $\overline{f}(0,r_k)=0$.  Thus, we let in the iteration step below, $p_0=r_k$.  We can see this with $f(z,w)=(97 z^7)+(100 z^2+77 z^3-77 z^4+64 z^5)w+(94-75 z^2)w^2+(-54)w^3$.  The first segment produces a characteristic equation of $97+100 x$.  When we process the function as above, $\overline{f}(z,w)=97+100w+zg(z,w)$ and therefore, the root of the characteristic equation is a root of $\overline{f}(0,w)$.
The iteration process is a variation of Newton iteration, a description of which can be found in Kung.  We let
\begin{equation}
\begin{aligned}
p_{n+1}&=p_n-\text{mod}\left(\frac{\overline{f}(z,p_n)}{d\overline{f}(p_n)},z^{2^{n+1}}\right), \\
p_0&=r_k,
\end{aligned}
\label{eqn035}
\end{equation}
where the mod function is simply extracting the $2^{n+1}-1$ terms of the Taylor series of the quotient.
After a set number of iterations we obtain the $n$-th iteration and write for the branch, 
\begin{equation}
\begin{aligned}
e_i&=\sum_{j=1}^i \lambda_i, \\
e_t&=\sum_{k=1}^t \lambda_k, \\
w_d(z)&=\frac{1}{z^{E}}\sum_{i=1}^{t} c_i z^{e_i}+z^{e_t} p_n(z^{1/d}),
\end{aligned}
\label{eqn036}
\end{equation}
where $E$ is the normal exponent of the function.
In this way, the method produces $n$ power series for an $n$-th degree function.
\end{enumerate}
\section{Implementing a numerical version of Newton Polygon in Mathematica}
Existing implementations of Newton polygon use exact arithmetic due to the possibility of numerical errors that produce incorrect polygon results.  However, these exact methods limit the types of functions that can be studied due to the computational complexity, storage space and execution times involved.  In this paper, a numeric version of the method is implemented and shown to produce acceptable results for a variety of different function types.  Numerical approximations however introduce the following issues:
\begin{enumerate}
\item{Multiple root resolution,}
\item{Residual coefficients close to zero,}
\item{Loss of accuracy due to coefficient expansion or contraction,}
\item{Numerical integration drift.}
\end{enumerate}
These issues are dealt with as follows:
\subsection{Multiple root resolution and residual coefficients}
Given the function:
$$f(z,w)=(w-1)(w-2)^2(w-3)^3(w-4)^4-z,$$
when we create the polygon for this equation, we obtain as the characteristic equation for $c_1$,
\begin{multline}
27648-110592 x+192384 x^2-192832 x^3+123852 x^4-53428 x^5\\
+15715 x^6-3118 x^7+400 x^8-30 x^9+x^{10}=0.
\end{multline}
Solving this for $x$ using machine-precision, the roots are 
$$\{1.,2.,2.,2.9993,3.00009,3.00103,3.9948,3.99966,4.00375,4.00589\}$$
so that if we used a tolerance of just $10^{-3}$, we would not pick up the multiple roots.  In order to successfully resolve the multiple roots, the algorithm must check the roots against the accuracy of the computation.  If the roots do not differ by this accuracy, then the algorithm used in this study considers the roots multiples.  Mathematica has features to monitor and control the accuracy and precision of all computations.  We therefore begin the algorithm by computing the singular points of the function with a high degree of precision.  This then determines the initial precision of succeeding calculations.  In this study, the precision of the singular points, $\{s_n\}$, was set between $400$ and $5000$ digits.  We then carefully monitor the accuracy of the computations.  If the number of significant digits or precision of a computation drops to zero, then the quantity is considered numerically zero.  Consider the two numbers:
$$
\begin{aligned}
\text{n1} &= 2.3333333333333333333345353\\
\text{n2} &= 2.3333333333333333333379889.
\end{aligned}
$$
These numbers differ in the 21'th decimal place.  Therefore, if the accuracy or number of significant digits to the right of the decimal place of the numbers is set to a value greater than 21, then the precision of the difference $n1-n2$ will be greater than zero.  However, if the accuracy of the numbers are both set to 20, then the precision of the difference $n1-n2$ drops to zero.  This is the principle used to check both for multiple roots and residual coefficients in the computations.  If the roots of the characteristic equation are computed to $150$ digits of accuracy and then that accuracy is manually decreased to $145$,  then multiple roots $r_1$ and $r_2$ when subtracted, will produce a number with zero significant figures or zero precision.  Likewise, if during the Newton polygon phase of the computation, the precision of the iterates, $f_i(z,w)$, drops to zero when the accuracy is manually decreased to an accuracy below its current accuracy, then that value is considered numerically zero and dropped.  This of course raises the possibility of numerical error if significant results are below this threshold.  However the algorithms designed in this paper were written with diagnostic tools to detect this and other possible errors due to numerical precision.  One tool is a table reporting the difference $|r_1-r_2|$ for all roots of the characteristic equation.  The user can visually inspect these values if there is a concern for numerical error.  Another tool reports all residual terms removed from the computations. These can also be checked and if necessary, the precision of the computations increased.  
\subsection{Accuracy reduction}
One limiting factor affecting the accuracy of the results is the computation of roots to the characteristic equation.  The accuracy of the roots is a function of the precision of the input data.  Another factor is the quotient $\displaystyle\frac{\overline{f}(z,p_n)}{d\overline{f}(p_n)}$ in the iteration phase of the calculation.  After each iteration, the accuracy of the results is usually less than the input data.  After several iterations, the accuracy can drop significantly.  In this implementation of the algorithm, the accuracy of all results are continuously checked and if it goes below zero, the algorithm stops and prompts the user to increase the precision of the calculations.

For complicated functions such as the twenty-degree function studied in Test \ref{testcase005}, the Newton iteration process produces extremely large coefficients on the order of several hundred.  Therefore, if the precision of the initial data was only $200$ digits, and the coefficients blow-up to values exceeding $10^{200}$, then the number of significant digits to the right of the decimal place drop to zero or even negative.  The solution to this problem is simply to increase the precision of the input data.  In order to detect when this occurs, the accuracy of the computations are checked and when the accuracy drops below zero, the user is notified to increase the working precision.
\subsection{Numerical integration limitations}
The function studied in Test \ref{testcase005} has a third-order pole located at the origin.  During the numerical integration phase of the continuation algorithm, the starting point for the integration is very close to the origin.  During the testing, the numerical results varied greatly from the expected results.  Although this problem may have been addressed by careful tuning of the numerical integrator, another remedy is to use the normalized function for the continuation analysis as was done in the test since the pole at the origin is removed.  Another potential problem is when the integration path encounters widely-varying function values or if the path is very long.  In Test \ref{testcase006}, the integration path is of the order of $10^{26}$.  If a numerical problem is suspected, the precision of the numerical integration can be increased in an effort to resolve it.  Another solution is to carefully monitor the integration as was done for Test \ref{testcase006} as described in the Conclusion section below.  
\section{Determining the radius of convergence of an algebraic power series}
\begin{figure}
	\centering
		\includegraphics[scale=0.75]{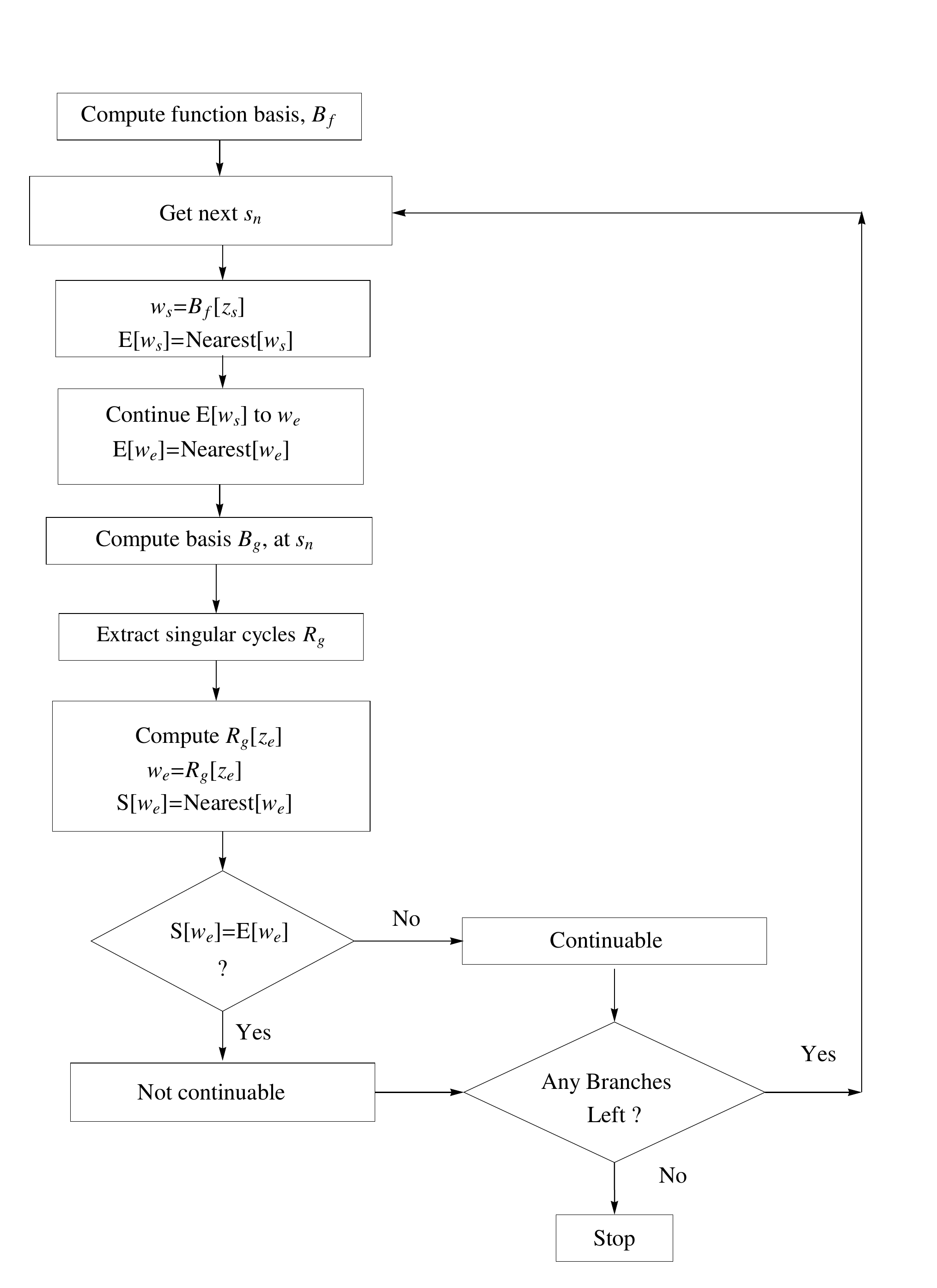}
	\caption{Branch Continuation Flowchart}
	\label{continuationFlowchart}
\end{figure}

\begin{figure}
	\centering
		\includegraphics[scale=0.5]{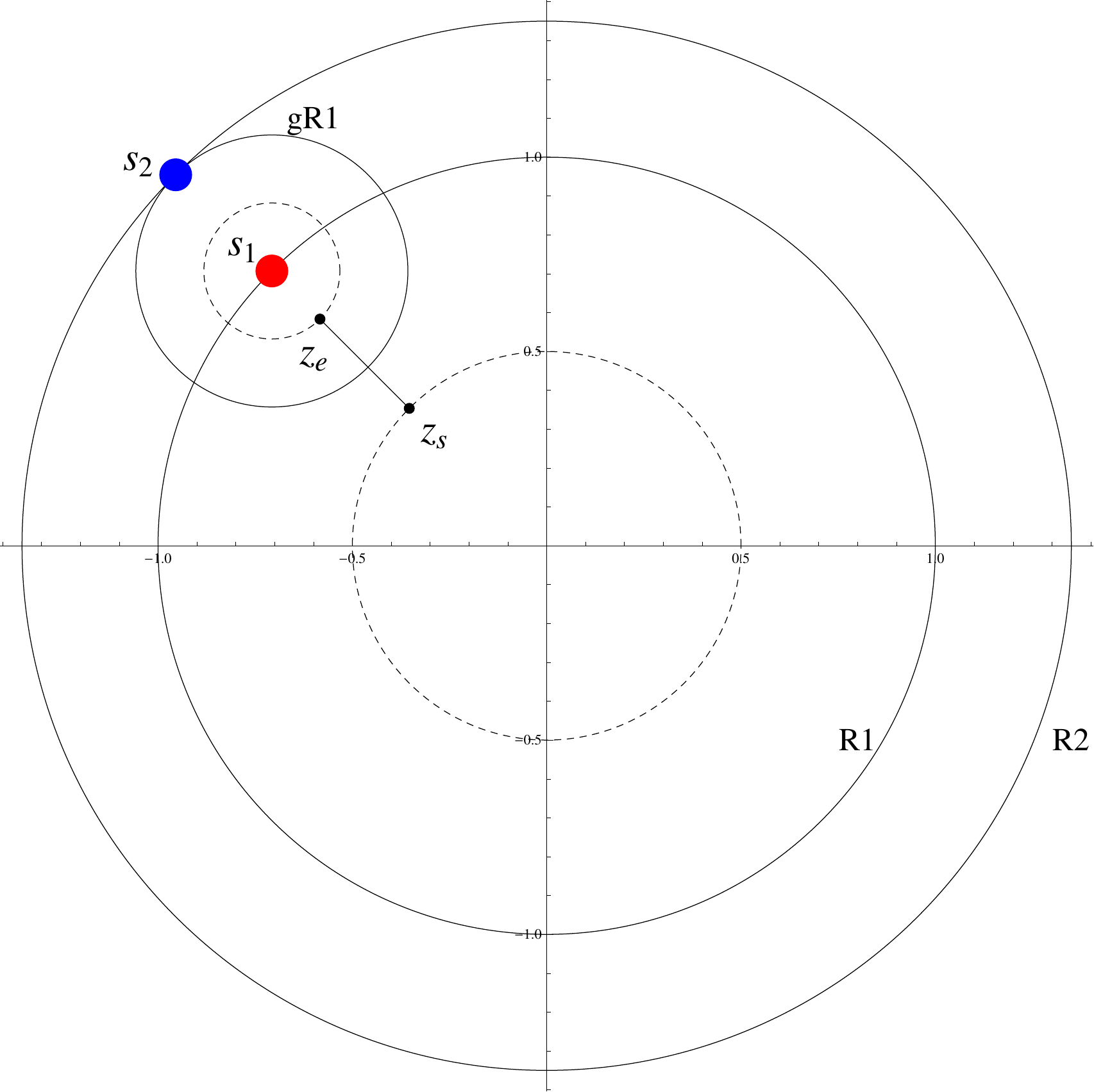}
	\caption{Branch Continuation Path}
	\label{figure025}
\end{figure}

The geometry of algebraic functions described above is used to determine the radii of convergence of their associated power series.  The process is depicted in  Figures $\ref{continuationFlowchart}$ and $\ref{figure025}$.  In general, the procedure checks each branch of the function for analytic continuity over singular points.  Figure $\ref{figure025}$ shows the setup over the first singular point in red which has a nearest neighbor in blue. 
The first step is to compute the expansion around the origin .  This expansion will converge at least to the nearest singular point with distance $R1$ in the figure.  The value of the function is then computed at one-half the distance of $R1$ at the point $z_s$ using the expansions computed above.  These values are then adjusted to the actual values as determined by the function $\textbf{Nearest}(f(z_s,w)=0)$.  The algorithm first determines the smallest difference between the roots of $f(z_s,w)=0$.  This value is given by $p_{\text{min}}$.  All comparisons are then made to a tolerance of $\displaystyle\frac{p_{\text{min}}}{N}$ where $N$ is an integer which in the test cases below, was $100$.  If this tolerance is not met at the current number of terms of the series, the the user is prompted to generate series with more terms.   The next step is to numerically extend each branch from $z_s$ in a radial direction to $z_e$.  This is done using the differential equation for the function with initial value at $z_s$.  These values are adjusted to the exact values at this point using the method described above.  The third step is to compute the local expansion around the selected singular point.  We know these expansions will converge in a radius at least equal to the distance to the nearest singular point or $gR1$.  Taking one-half the distance, we compute the value of each local branch at $z_e$ and adjust for exact values.  If the analytically-continued values from $z_s$ to $z_e$ impinge upon ramified branches or poles at $s_n$, and if the singular points were processed in sequential order, the power series has a radius of convergence equal to this singular ring.  If a branch at the origin does not impinge onto local singular branches for all singular points in a ring, the power series for this branch has a radius of convergence beyond this ring.  We continue in this fashion until no further continuations are possible.  The details of each step are as follows:\\
\begin{enumerate}
\item{Computing the singular points $\{s_n\}$:}

 In this study, the singular points were computed to sufficient precision to mimimize the potential for possible numerical errors described above.  This precision was usually between $400$ and $5000$.
\item{Computing the function basis $B_f$:}

The function basis is the set of conjugately-distinct power series centered at the origin.  Thus for a $4$-cycle, the basis would have only one series representing this branch.  Likewise for the other multi-valued branches.  As many as $4035$ terms of a series were generated in this study.
\item{Computing $B_f(z_s)$ and $E(w_s)$:}

A line is drawn from the origin to the selected singular point.  On this line and one-half the distance to the first singular point, the point $z_s$ is chosen.  Using the basis computed above, we compute the values of the function branches at $z_s$.  These values will differ from the exact values (to numerical precision) of $f(z_s,w)=0$, depending on how many terms of the series are used.  The values are adjusted by comparing the computed value $w_d(z_s)$ to the ideal values and choosing the closest match and then checking that the computed values are less than $p_{\text{min}}/100$.  The exact values are $E(w_s)$. 
\item{Continue $E(w_s)$ to the point $z_e$:}

For each singular point, the distance to it's nearest neighbor is computed.  We choose a point which is one-half this distance and on the line to the origin and mark a point $z_e$.  If this point is smaller in absolute value to $z_s$, it is adjusted (made smaller) so that $|z_e|-|z_s|\geq 0$.  Now using the differential equation for the function,
$$
\frac{dw}{dr}=-\left(\frac{f_z}{f_w}\right) \frac{dz}{dr},
$$  
we numerically continue each branch from $z_s$ to the point $z_e$, and as was done with $z_s$, adjust the values to the closest match, $E(w_e)$.  
\item{Compute basis $B_g$ at selected singular point:}

We next compute the basis around the chosen singular point $s_n$ by making the substitution $g(z,w)=f(z+s_n,w)$ and thus an expansion around zero for $g(z,w)$ is an expansion of $f(z,w)$ around $s_n$. The singular points are chosen in sequential order according to their radial distance from the origin in order to avoid integrating over a singularity.  This could occur if two singular points are in alignment with the numerical integration path and analyticity over the closest singular point was not checked prior to checking the more distant point.  The basis, $B_g$,  will contain at least one singular branch.  This singularity will either be a pole and or a ramified branch with multiple values.  Poles and ramified branches are barriers to extending the radius of convergence of power expansions across singular points and therefore if a  branch from the basis $B_f$ impinges upon any of the singular branches of $B_g$ at $s_n$, its power series has a ring of convergence equal to the ring number of the current ring, $r_n$.  If a branch does not impinge upon the singular branches $B_g$ for all singular points of a ring, then the series has a radius of convergence which extends across this singular point.  We therefore select the singular branches from $B_g$, compute $B_g(z_e)$ for each sheet of the branch, adjust those values to their exact values, and then compare those values to each of the continued values $B_f(z_e)$.  Any branches in which $B_f(z_e)=B_g(z_e)$ to numerical precision has a power series with a radius of convergence equal to $|s_n|$.  Branches which do not agree are then checked against the next singular point.
\item{Continue the above process until no remaining series can be extended over the next singular point.}

\end{enumerate}
\section{Algorithm Testing}

In order to obtain some empirical measure of the accuracy of the test results, all power series were checked with values $z_c$ and $z_d$ such that $|r_{n-1}|<|z_c|<|r_n|$ and $|r_n|<|z_d|<|r_{n+1}|$ with adjustments made if $r_c$ is the first or last ring.  If the computed ring of convergence is correct and the series are precise, the partial sums will exhibit convergent behavior in the former case and divergent behavior in the later case.  Additionally, the series were checked against $100$ random points in the range $\frac{|r_c|}{100}<|z_i|<\frac{99}{100}|r_c|$ with the values compared against the expected value given by the function $\textbf{Nearest}(f(z_i,w)=0)$.  The maximum error of the set is reported in the tables. 
\subsection{Test Cases}
\begin{enumerate}
\item{Function with cycles one through four:}
\label{testcase001}
\begin{equation}
\begin{aligned}
f(z,w)&=(z^{14}+3 z^{15})+(2 z^{15}+2 z^{16})w+(3 z^{15}-20 z^{16})w^2+(4 z^{15})w^3\\
&+(10 z^5-z^6+2 z^7)w^4+(8 z^{10})w^5+(9 z^{10})w^6\\
&+(20 z)w^7+(3 z)w^8+(2)w^9+(5)w^{10}.
\end{aligned}
\label{eqntestcase001}
\end{equation}
Continuation results for function $(\ref{eqntestcase001})$ are shown in Table $\ref{table001}$.  The test was run with a precision of $800$ using $54$ terms of each series.  Note the continuation column.  In the first row, the four and three cycle branches are continuable over the first singular point.  These series then have a radius of convergence at least equal to the absolute value of the second singular point or approximately 0.3329.  Note however, only the $3$-cycle branch is continuable across the singular points in the second and third rings but not the fourth.  This means $w_3$ has a radius of convergence equal to the size of the fourth ring or approximately 0.636.

Convergence results are tabulated in Table $\ref{table002}$ with $r_c$ representing the ring of convergence as described above.   The partial sum analysis was run for each branch and agreed with the computed radii of convergences in the table.  Divergent behavior for $w_3$ and $w_4$ was not observed until both partial sums had accumulated more than $2000$ terms.
\begin{center}
\begin{table}
\caption{\text{Branch Continuations for Test Case} $\ref{testcase001}$}
\label{table001}
\begin{tabular}{|c|c|c|c|}
\hline
 \text{Ring} &  \text{Singularity} & \text{Abs Value} & \text{Continuations} \\
 \hline
 1 &   -0.002469 & 0.002469 & $\{w_4,w_3\}$ \\
 2 &  -0.3329+0.00083i & 0.3329 & $\{w_3\}$ \\
   &  -0.3329-0.00083i & 0.3329  & $\{w_3\}$ \\
 3 &  -0.3341          & 0.3341  & $\{w_3\}$ \\
 4 &  0.244+0.587 i    & 0.636   & $\{\}$\\
\hline
\end{tabular}
\vspace{10pt}
\end{table}
\end{center}
%
% ROC results 
%
\begin{center}
\begin{table}
\caption{\text{Radius of Convergence Results for Test $\ref{testcase001}$}}
\label{table002}
\begin{tabular}{|c|c|c|c|c|}
\hline
\text{Cycle} & $r_c$  & $|r_c|$ & \text{Terms} & \text{Max error}\\
\hline 
 $w_1$ &  1 & $0.00247$ & $256$ & $10^{-7}$\\
 $w_2$ &  1 & $0.00247$ & $254$ & $10^{-6}$\\
 $w_3$ &  4 & $0.6363$ & $4094$ & $10^{-14}$\\
 $w_4$ &  2 & $0.3329$ & $4086$ & $10^{-19}$ \\
    \hline
\end{tabular}
\vspace{10pt}
\end{table}
\end{center}%
\item{Function having cycles 1, 2, 3, 4, and 5:}
\label{testcase002}
\begin{equation}
f(z,w)=(z^{30}+z^{32})+(z^{14}+z^{20})w^5+(z^5+z^9)w^9+(z+z^3)w^{12}+(6)w^{14}+(2+z^2)w^{15}.
\label{eqn037}
\end{equation}
Results of this test are in table \ref{tabletestcase002}.  This function has a $1$-cycle branch which extends across $61$ singular rings containing $118$ singular points.  Cycles one through three were tested with $1024$ terms at $1000$ digits of precision.  In order to obtain definitive convergence/divergence behavior of the partial sum plots for the $5$-cycle branch, $4035$ terms were used because the partial sums initially exhibited convergent behavior but began diverging after approximately $2500$ terms.  
\begin{center}
\begin{table}
\caption{\text{Radius of Convergence Results for Test $\ref{testcase002}$}}
\label{tabletestcase002}
\begin{tabular}{|c|c|c|c|c|}
\hline
\text{Cycle} & $r_c$  & $|r_c|$ & \text{Terms} & \text{Max error}\\
\hline 
 $w_1$ &  61 & $1.093$ & $1024$ & $10^{-15}$\\
 $w_2$ &  1 & $0.1168$ & $1024$ & $10^{-12}$\\
 $w_3$ &  1 & $0.1168$ & $1024$ & $10^{-10}$\\
 $w_4$ &  4 & $0.505$ & $998$ & $10^{-7}$ \\
 $w_5$ &  14 & $0.6413$ & $4035$ & $10^{-10}$\\
   \hline
\end{tabular}
\vspace{10pt}
\end{table}
\end{center}
\item{Multiple $n$-cycle branches:}
\label{testcase003}
\begin{equation}
f_2(z,w)=z w^4\left[a(1-\overline{a}^2 w^2)\right]^4-\left[\overline{a}(a^2-w^2)\right]^4,\quad a=3-14i.
\label{eqn031}
\end{equation}
This function has three $4$-cycle branches which cannot be extended past the first singular point. This implies that the first singular point must have at least three ramified sheets and in fact the function ramifies into two $2$-cycle branches (with additional single-cycle branches) at this singular point.  Convergence results are in Table \ref{tabletestcase003}.  Branch $w_{4,3}$ is an example of a ramified pole.  Partial sum testing confirmed the computed rings of convergence.  

The maximum error observed for these series is relatively large.  This can be explained due to the fact that the radius of convergence is extremely small, on the order of $10^{-16}$.  The random points used to check the series were therefore all very close to the limit of convergence for each series and thus a relatively large difference between the series value and actual value would be expected for the small number of terms tested.
\begin{center}
\begin{table}
\caption{\text{Radius of Convergence Results for Test $\ref{testcase003}$}}
\label{tabletestcase003}
\begin{tabular}{|c|c|c|c|c|}
\hline
\text{Cycle} & $r_c$  & $|r_c|$ & \text{Terms} & \text{Max error}\\
\hline 
 $w_{4,1}$ &  $1$ &  $2.96\times 10^{-16}$ & 257 & $10^{-4}$\\
 $w_{4,2}$ &  $1$ &  $2.96\times 10^{-16}$ & 257 & $10^{-4}$\\
 $w_{4,3}$ &  $1$ &  $2.96\times 10^{-16}$ & 128 & $10^{-4}$\\
   \hline
\end{tabular}
\vspace{10pt}
\end{table}
\end{center}
\item{Function with a recursive Newton polygon:}
\label{testcase004}
\begin{equation}
f(z,w)=((w^3 + z^2)^2 + z^3 w^2)^2 + z^7 w^3.
\label{eq031}
\end{equation}
This function has two $6$-cycle branches and when analyzed with Newton polygon, produces a 3-level recursive polygon tree which means it generates multiple roots for two polygon phases.    This checks the algorithm's ability to detect multiple roots at a level other than the first when the coefficients are exact with infinite precision.  Results of this test are shown in Table $\ref{tabletestcase004}$ using $1026$ terms of each series.
\begin{center}
\begin{table}
\caption{\text{Radius of Convergence Results for Test $\ref{testcase004}$}}
\label{tabletestcase004}
\begin{tabular}{|c|c|c|c|c|}
\hline
\text{Cycle} & $r_c$  & $|r_c|$ & \text{Terms} & \text{Max error}\\
\hline 
 $w_{6,1}$ &  $1$ &  $0.9585$ & 1026 & $10^{-10}$\\
 $w_{6,2}$ &  $1$ &  $0.958$ & 1026 & $10^{-12}$\\
   \hline
\end{tabular}
\vspace{10pt}
\end{table}
\end{center}
\item{$20$-degree function with third-order pole at origin}
\label{testcase005}
\begin{equation}
\begin{aligned}
f(z,w)&=(-14 z-68 z^4+83 z^5+88 z^6)\\
&+(-20+19 z-19 z^2+42 z^6+25 z^7+54 z^8+66 z^9)w\\
&+(-16 z+59 z^2+54 z^3+36 z^5-91 z^6-14 z^7)w^2\\
&+(47-17 z^5+64 z^6+94 z^8+68 z^9)w^3\\
&+(-62 z+27 z^2-25 z^3+39 z^4)w^4\\
&+(-4+2 z+11 z^5-13 z^6+85 z^{10})w^5\\
&+(53 z+5 z^2-65 z^3+57 z^7-75 z^{10})w^6\\
&+(-13+92 z+23 z^2+z^4-15 z^5+23 z^7)w^7\\
&+(85 z-39 z^4+78 z^5+48 z^7-26 z^9+2 z^{10})w^8\\
&+(2 z^2-56 z^3+9 z^{10})w^9\\
&+(-17 z^2-65 z^3+77 z^4+64 z^8-45 z^9+96 z^{10})w^{10}\\
&+(39+84 z^3+90 z^4-6 z^5+6 z^6-57 z^8+39 z^9)w^{11}\\
&+(-48 z+79 z^2-22 z^5+75 z^6-3 z^9)w^{12}\\
&+(55+100 z^2-58 z^6-19 z^7+83 z^9-41 z^{10})w^{13}\\
&+(67+85 z-8 z^3+16 z^7)w^{14}\\
&+(-23 z^5-39 z^7-20 z^{10})w^{15}\\
&+(56+76 z+57 z^2+100 z^3-40 z^4+68 z^5-55 z^7+50 z^9+52 z^{10})w^{16}\\
&+(69-42 z+53 z^2-89 z^4-13 z^7+55 z^8)w^{17}\\
&+(-12 z^2-22 z^8-11 z^{10})w^{18}\\
&+(62+23 z-4 z^2-99 z^4+9 z^5-99 z^6+57 z^9-90 z^{10})w^{19}\\
&+(92 z^3-91 z^8+63 z^9)w^{20}.
\end{aligned} 
\label{eqn032}
\end{equation}
This function exemplifies what can happen when numerically integrating near a pole.  In step $4$ of the continuation process, we numerically integrate from $z_s$ to $z_e$.  However, if $z_s$ is near a pole the numerical integration may suffer, as in this case, due to the large derivative involved.  Although it may be possible to minimize this problem with efficient use of the numerical integrator, we can avoid the problem by continuation over the normalized function.  Recall, the normalization process removes the singular point at the origin.  When normalized though, the powers on $z$ are usually increased.  In this particular function, the largest power of $z$ rises to $63$.  This causes the coefficients of the iteration functions, $p_n$, to grow extremely large and thus, to maintain accuracy of the data, the working precision must be increased.  In this particular case, it was increased to $1000$. Table \ref{tabletestcase006} summarizes the results of this test.   The partial sum study for this function agreed with all convergence results.
\begin{center}
\begin{table}
\caption{\text{Convergence results for Test $\ref{testcase005}$}}
\label{table005}
\begin{tabular}{|c|c|c|c|c|}
\hline
 \text{Cycle} & $r_c$  & $|r_c|$ & \text{Terms} & \text{Max Error}\\
 \hline
 $w_{1,11}$ & 3 & 0.249329 & 256 & $10^{-6}$ \\
 $w_{1,2}$ & 2 & 0.139071 & 256 & $10^{-8}$\\
 $w_{1,3}$ & 2 & 0.139071 & 256 & $10^{-7}$\\
 $w_{1,4}$ & 4 & 0.251398 & 256 & $10^{-11}$\\
 $w_{1,5}$ & 4 & 0.251398 & 256 & $10^{-9}$\\
 $w_{1,6}$ & 2 & 0.139071 & 256 & $10^{-10}$\\
 $w_{1,7}$ & 17 & 0.578850 & 256 & $10^{-7}$\\
 $w_{1,8}$ & 17 & 0.578850 & 256 & $10^{-7}$\\
 $w_{1,9}$ & 1 & 0.0430598 & 256 & $10^{-12}$\\
 $w_{1,10}$ & 1 & 0.0430598 & 256 & $10^{-7}$\\
 $w_{1,12}$ & 1 & 0.0430598 & 256 & $10^{-8}$\\
 $w_{1,13}$ & 1 & 0.0430598 & 256 & $10^{-10}$\\
 $w_{1,14}$ & 6 & 0.267128 & 256& $10^{-8}$\\
 $w_{1,15}$ & 6 & 0.267128 & 256 & $10^{-8}$\\
 $w_{1,16}$ & 5 & 0.257919 & 256 & $10^{-7}$\\
 $w_{1,17}$ & 5 & 0.257919 & 256 & $10^{-6}$\\
 $w_{1,18}$ & 5 & 0.257919 & 256 & $10^{-7}$\\
 $w_{1,19}$ & 6 & 0.267128 & 256 & $10^{-12}$\\
 $w_{1,20}$ & 6 & 0.267128 & 256 & $10^{-7}$\\
 $w_{1,1}$ & 13 & 0.530749 & 256 & $10^{-5}$\\
  \hline
\end{tabular}
\vspace{10pt}
\end{table}
\end{center}
\item{Function of degree $55$:}
\label{testcase006}
\begin{equation}
f(z,w)=z-\prod_{j=1}^{10} (w-j)^j.
\label{eqn050}
\end{equation}
This function represents an inverse of a $55$-degree polynomial and was analyzed with $128$ terms at a working precision of $400$.  Also, the working precision and accuracy of the numerical integration was set to $40$ and $30$.  Results of this test case are shown in Table $\ref{tabletestcase006}$ and agreed with the partial sum testing.
\begin{center}
\begin{table}
\caption{\text{Convergence Results for Test $\ref{testcase006}$}}
\label{tabletestcase006}
\centering
\begin{tabular}{|c|c|c|c|c|}
\hline
\text{Cycle} & $r_c$  & $|r_c|$ & \text{Terms} & \text{Max Error}\\
\hline 
 $w_1$ & 9 & $2.39\times 10^{38}$ & $256$ & $10^{-12}$  \\
 $w_2$ & 8 & $2.9\times 10^{32}$ & $256$ & $10^{-6}$\\
 $w_3$ & 7 & $4.79 \times 10^{26}$ & $256$ & $10^{-7}$\\
 $w_4$ & 6 & $1.95\times 10^{21}$ & $256$ & $10^{-6}$\\
 $w_5$ & 5 & $ 3.75\times 10^{16}$ & $256$ & $10^{-6}$\\
 $w_6$ & 3 & $7.24\times 10^{12}$ & $256$ & $10^{-6}$ \\
 $w_7$ & 2 & $3.69\times 10^{10}$ & $256$ & $10^{-6}$\\
 $w_8$ & 1 & $2.12\times 10^{10}$ & $256$ & $10^{-6}$\\
 $w_9$ & 1 & $2.12\times 10^{10}$ & $256$ & $10^{-5}$\\
 $w_{10}$ & 4 & $2.34\times 10^{13}$ & $256$ & $10^{-4}$\\
  \hline
\end{tabular}
\vspace{10pt}
\end{table}
\end{center}
\item{Function with continuations across single-sheet poles}:
\label{testcase007}
\begin{equation}
\begin{aligned}
f(z,w)&=(3+4 z)+(-6 z^2-\frac{3 z^5}{2})w+(\frac{1}{2}-16 z+\frac{7 z^2}{8})w^2+(\frac{3}{4}-2 z+12 z^4)w^3\\
&+(15+\frac{z^2}{3}+\frac{22 z^3}{15})w^7+(-\frac{1}{1000}-\frac{z}{25}+\frac{z^2}{2}-\frac{z^3}{5}+2 z^4)w^8.
\end{aligned}
\end{equation}
The continuation table of this function is given in Table \ref{tabletestcase007}.  The poles are in red and since all the branches are single-sheets, the notation $w_{d,n}$ is dropped and only the sort orders are reported.  All but $w_{1,8}$ are continuable across the first two poles and branch $w_{1,1}$ continues over all  poles of the function.  Convergence results are in Table \ref{tabletestcase007b}.  Cycle $w_{1,8}$ at $256$ terms did not produce a low maximum error at $256$ terms.  This may have been due to the fact that this branch sheet is $15000$ at the origin and quickly reaches $140,000$ at the boundary of the testing data.  However, at $2048$ terms, the maximum error was reduced to $10^{-9}$.
\begin{center}
\begin{table}
\caption{\text{Continuation results for Test $\ref{testcase007}$}}
\label{tabletestcase007}
\centering
\begin{tabular}{|c|c|c|c|}
\hline
\text{Ring} & $s_i$  & $|s_i|$ & \text{Continuations}\\
\hline 
1 & $\textcolor{red}{-0.019968}$ & $0.019968$ & $\{1,2,3,4,5,6,7\}$ \\
2 & $\textcolor{red}{0.1}$ & $0.1$ & $\{1,2,3,4,5,6,7\}$ \\
3 & $-0.290-0.155i$ & $0.329$ & $\{1,3,5,6,7\}$ \\
  & $-0.290+0.155i$ & $0.329$ & $\{1,3,5,6,7\}$ \\
4 & $-0.353-0.346i$ & $0.494$ & $\{1,6\}$\\
  & $-0.353+0.346i$ & $0.494$ & $\{1\}$ \\
5 & $\textcolor{red}{0.1-0.5i}$ & $0.5004$ & $\{1\}$ \\
  & $\textcolor{red}{0.1+0.5i}$ & $0.5004$ & $\{1\}$ \\
6 & $-0.0733-0.543i$ & $0.548$ & $\{ \}$ \\
  \hline
\end{tabular}
\vspace{10pt}
\end{table}
\end{center}
\begin{center}
\begin{table}
\caption{\text{Convergence Results for Test $\ref{testcase007}$}}
\label{tabletestcase007b}
\centering
\begin{tabular}{|c|c|c|c|c|}
\hline
\text{Cycle} & $r_c$  & $|r_c|$ & \text{Terms} & \text{Max Error}\\
\hline 
 $w_{1,1}$ & 6 & $0.5489$ & $256$ & $10^{-7}$  \\
 $w_{1,2}$ & 3 & $0.3288$ & $256$ & $10^{-9}$\\
 $w_{1,3}$ & 3 & $0.3288$ & $256$ & $10^{-8}$\\
 $w_{1,4}$ & 3 & $0.3288$ & $256$ & $10^{-7}$\\
 $w_{1,5}$ & 3 & $0.3288$ & $256$ & $10^{-7}$\\
 $w_{1,6}$ & 4 & $0.4943$ & $256$ & $10^{-7}$ \\
 $w_{1,7}$ & 4 & $0.4943$ & $256$ & $10^{-7}$\\
 $w_{1,8}$ & 1 & $0.01997$ & $2048$ & $10^{-9}$\\
  \hline
\end{tabular}
\vspace{10pt}
\end{table}
\end{center}
\end{enumerate}
\section{Timing Statistics}
The algorithms were run on a dual-Pentium processor machine at 2.2 GHz under Mathematica version 8.0.
Table $\ref{tab:table0025}$ summarizes typical timing results for the test cases. S is the timing to compute singular points and some auxiliary tables, $N_P$, the timing for the Newton polygon algorithm, and $A_C$, the timing for the  continuation algorithm.
 \begin{table}
    \caption{Run-time Statistics (in seconds)}
    \label{tab:table0025}
    \centering
    \begin{threeparttable}
      \begin{tabular}{
        cccccc
              }\toprule
        {\multirow{2}{*}[-0.5ex]{Test Case}} & {\multirow{2}{*}[-0.5ex]{$p_w$}} & {\multirow{2}{*}[-0.5ex]{Terms}} & \multicolumn{3}{c}{Timing} \\ \cmidrule(lr){4-6}
        & & & S & {$N_P$} & {$A_C$} \\ \midrule
        1 & 800 & 54 & 1.5 & 0.17 & 21.5 \\
        2 & 800 & 17 & 11.7 & 0.17 & 4605 \\
        3 & 2500 & 32 & 0.1 & 40.2 & 21.4 \\
        4 & 800 & 66 & 0.1 & 9.3 & 6.5\\
 $*5$ & 1000 & 63 & 150.2 & 0.8 & 2871\\
 6 & 400 & 128 & 0.5 & 793 & 4181 \\ 
 7 & 1500 & 256 & 0.8 & 63 & 155\\ \bottomrule
       \end{tabular}
      \begin{tablenotes}
        \footnotesize
        \item[$*$] Results for normalized function
      \end{tablenotes}
    \end{threeparttable}
  \end{table}
\begin{center}
\begin{table}
\caption{\text{NDSolve results for Test $\ref{testcase006}$, singular point $7$}}
\label{conclusiontable001}
\centering
\begin{tabular}{|c|c|c|c|c|}
\hline
\text{Cycle} &\text{Sheet} & \text{Actual Value}  & \text{NDSolve results} & \text{Difference}\\
\hline 
 $w_1$ & $1$ & $1.000000$ & $1.000000$ & $0.\times 10^{-40}$  \\
 $w_2$ & $1$ & $2.000000-0.00005i$ & $2.000000-0.00005i$ & $8.79\times 10^{-20}$ \\
       & $2$ & $2.000000+0.00005i$ & $2.000000+0.00005i$ & $8.79\times 10^{-20}$ \\
 $w_3$ & $1$ & $2.9637-0.0433i$  & $2.9637-0.0433i$ & $1.137\times 10^{-19}$\\
       & $2$ & $3.1006$ & $3.1006$ & $1.555\times 10^{-19}$\\
       & $3$ & $2.9637+0.0433i$ & $2.9637+0.0433i$ & $1.137\times 10^{-19}$\\
  \hline
\end{tabular}
\vspace{10pt}
\end{table}
\end{center}
\section{Conclusions}
The limiting factor controlling radii of convergence of algebraic power series is branch-sheet continuation to either a ramified covering or single-sheet pole.  This was the principle factor used in the continuation algorithm described in this paper.  For each test case studied, the partial-sum testing agreed with the computed radii of convergence for all branches: the partial sums exhibited convergent behavior for $|r_{c-1}|<|z_c|<|r_c|$, and showed divergent behavior in for $|r_c|<|z_d|<|r_{c+1}|$.  And since the singular points are the roots  to the resultant polynomial with rational coefficients, these values can be computed to arbitrary precision for a wide class of functions.

However, the software may require tuning for a particular function.  For example, the algorithms rely on determining an exact match for function values when these values are only computed approximately.  The matching is done by finding the closest match to within a tolerance of $\frac{1}{100}$ of the smallest separation between function values.  The code will detect when this tolerance is not met and will prompt the user to either generate more terms of the series in the Newton polygon phase or to numerically integrate with a higher working precision.  Other potential problems could come from extremely close polynomial roots, very small or large coefficients encountered during the analysis, algebraic branches which have their sheets extremely close to one another, and other errors associated with numerical integration of a differential equation as was shown in Test $\ref{testcase005}$.  Steps to identify these errors were described above and diagnostic tools were designed to identify possible problems. 

Consider Test \ref{testcase006}.  The singular points are extremely large and therefore, numerical integration from the point $z_s$ to $z_e$ has to traverse increasingly large distances.  One branch extends out to  $10^{38}$.  How can we be certain the integration does not over-track the branch and land on a different sheet along the integration path?  One way to minimize this is to integrate over a well-behaved section of the function.  This function has no poles and so is well-behaved along the integration paths. Another way is to integrate with a high degree of accuracy.  The test case used an accuracy goal of $30$ digits.  And still a third way of handling this possible problem is with diagnostic tools.  The code written for this work has an option for reporting the results of the integration.  Table $\ref{conclusiontable001}$ is an example of one such diagnostic tool.  It was generated during the analysis of Test \ref{testcase006} about the seventh  singular point near $-4.97\times 10^{26}$.  The values are however reported to six decimal places to avoid clutter and so some numeric quantities for this particular function appear to either be identical or zero but are actually different at a greater precision.  The important point is the last column which reports the difference between the actual value of the function at $z_e$ and the value determined by numerical integration.  The table reports that after integrating a distance of approximately $10^{16}$, the integration did not vary by more than $10^{-19}$ and since the maximum distance between the sheets at $z_e$ was approximately unity, this small difference compared to the sheet separations lends credence to the accuracy of the results.  In this case, the sheets of the function were basically flat.  However, had the sheets values varied greatly in the integration intervals, the analysis may have encountered problems.

However, as stated earlier, one objective of this work was to produce a first version of a software tool that could successfully compute Puiseux series and their radii of convergence for a wide variety of functions not accessible by exact arithmetic means.  As best as could be determined, the algorithm produced good results for the test cases that were studied.  If all phases of the analysis are carefully monitored with the available diagnostic tools, the potential for errors can be reduced.  Still though, the code has much room for improvement and can provide a means of computing radii of convergence precisely for the types of functions studied in this paper.  A greater variety of function types such as those with more polygon recursions, functions with higher ramified singular points, more complicated coefficients, or other morphologies to stress-test the software would produce a more robust algorithm as would a more careful analysis of the accuracy of the computations.  

One might consider using the partial sum testing to conceivably compute these convergence radii.  Simply compute a highly accurate series with many terms and check it's partial sum behavior in successive rings until it begins to diverge.  However, the point at which the series begins to diverge is not known and as was observed in the test, the partial sum may even appear to slowly converge and then start to diverge.  This is what happened for the four and three-cycle branches in Test \ref{testcase001}. 

In principle, we should be able to dispense with the numerical integration step of the continuation algorithm since the convergence domain of the power series being checked, and convergence domain of the power series centered at the chosen singular point will always overlap using the method described in this paper.  We would then need only compare the values from two different power series.  Doing this would eliminate the numerical integration step in the algorithm thus avoiding potential numerical errors.  However, this would entail the possibility of computing the value of a power series at a point very close to it's radius of convergence and near a singularity so might require the need for many terms to achieve a desired accuracy.

\end{document}